\documentclass[11pt]{article}
\usepackage{epsfig}
\usepackage{amssymb}
\usepackage{amscd}
\usepackage{epsf}
\usepackage{amsfonts,amssymb,amsthm,amsmath}

\long\def\comment#1\endcomment{\relax}

\newcounter{subsubsubsection}
\newcounter{subsubsubsubsection}

\makeatletter
\@addtoreset{subsubsubsection}{subsubsection}
\@addtoreset{subsubsubsubsection}{subsubsubsection}

\makeatother

\newcommand{\sevafigc}[4]{\begin{figure}[h]\centerline{
 \epsfig{file=#1,width=#2,angle=#3}}
\bigskip\caption{#4}\end{figure}}

\DeclareMathSizes{11.1}{10}{8}{6}

\DeclareMathOperator{\Hom}{Hom}

\newtheorem*{theorem*}{Theorem}

\newtheorem*{lemma}{Lemma}

\newtheorem*{proposition}{Proposition}

\newtheorem*{corollary}{Corollary}

\theoremstyle{remark}

\newtheorem*{remark}{Remark}

\newtheorem*{example}{Example}

\newtheorem*{question}{Question}

\theoremstyle{definition}

\newtheorem*{definition}{Definition}

\DeclareMathOperator{\Id}{Id}

\DeclareMathOperator{\Alt}{Alt}

\newcommand{\mb}{{\bullet}}

\newcommand\End{\mathrm{End}}

\newcommand{\Conf}{\mathrm{Conf}}

\newcommand{\g}{\mathfrak{g}}

\newcommand{\K}{\mathrm{K}}
\newcommand{\X}{\mathrm{H}}
\newcommand{\Der}{\mathrm{Der}}

\newcommand{\T}{\mathrm{T}}
\newcommand{\St}{\mathrm{St}}

\def\wtilde#1{\widetilde{#1}\vphantom{#1}}

\title{{\tt {\huge The CROCs, non-commutative deformations, and (co)associative bialgebras}}}
\author{{\tt {\LARGE Boris Shoikhet}}}

\begin{document}\maketitle

\begin{abstract}
{\tt We compactify the spaces $K(m,n)$ introduced by Maxim Kontsevich. The initial idea
was to construct an $L_\infty$ algebra governing the deformations of a (co)associative
bialgebra. However, this compactification leads not to a resolution of the
PROP of (co)associative bialgebras, but to a new algebraic structure we call here a CROC.
It turns out that these constructions are related to the non-commutative
deformations of (co)associative bialgebras. We construct  an associative dg algebra
conjecturally governing the non-commutative deformations of a bialgebra.
Then, using the Quillen duality, we construct a dg Lie algebra conjecturally
governing the commutative (usual) deformations of a (co)associative
bialgebra.

Philosophically, the main point is that for the associative bialgebras
the non-commutative deformations is maybe a more fundamental object than the usual
commutative ones. }
\end{abstract}

\section*{Introduction}
The non-commutative deformations of an object, although are not rigorously
defined at the moment, exist from a more philosophical point of view.
We suppose that the whole formal neighborhood of the
moduli space is some formal non-commutative manifold $X$, while its
commutative part is a subspace $X_0$. It is supposed that in a smooth point,
the algebra $A$ of functions on $X$ is a free associative algebra, and the
imbedding $X_0\hookrightarrow X$ gives the corresponding map of algebras
$p\colon A\to A_0$ where $A_0$ is the commutative algebra of functions on $X_0$.
Moreover, we suppose that the map $p$ is the projection
$p\colon A\to A/[A,A]=A_0$.

The deformations of an object are described by a dg Lie algebra $\g^\mb$ of
derivations of an appropriate resolution of the object. On the other hand,
this dg Lie algebra is related to the algebra $A_0$ in the following way:
the 0-th Lie algebra cohomology $H^0(\g^\mb)\simeq A_0$ (we do not consider
here the questions of completions, etc.). Very informally, we can say that
the dg Lie algebra $\g^\mb$ is the Quillen dual to the commutative algebra
$A_0$ (in reality, it is not true, it is true only for $H^0$). Keeping
this point of view, it is natural to describe the non-commutative
deformations by an object Quillen dual to $A$. The Quillen duality maps dg
commutative algebras to dg Lie algebras and vise versa, and it maps
associative dg algebras to itself (this duality is also known as the Koszul
duality for operads). This duality associates to an associative algebra $A$
the tensor algebra $\T(A^*[1])$ with the differential $\delta\colon A^*\to
(A^*)^{\otimes 2}$ which is dual to the product in $A$. Then, applying the
Quillen construction again, we obtain an associative dg algebra which 0th
cohomology is isomorphic to the initial algebra $A_0$.

On the other hand, the deformation dg Lie algebra $\g^\mb$ can be
constructed differently from the OPERADic point of view. The well-known
example is the deformation theory of associative algebras. The associative
algebras itself can be described as algebras over an operad
$\mathrm{Assoc}$. Then there is a geometric construction of a free minimal
resolution of the operad $\mathrm{Assoc}$. Namely, denote by $\St_n$ the
$n$-th Stasheff associahedron (of dimension $n-2$). Then this minimal model
$M^\mb$
is the direct sum $M^\mb=\bigoplus_{n\ge 2}C_\mb(\St_n)$ with the chain
differential where $C_\mb$ denotes the chain complex with the Stasheff cell
decomposition. By definition, an $A_\infty$ algebra is an algebra over this
free operad. In other words, for a vector space $V$, a map of operads
$M^\mb\to\End(V)$ is the same that an $A_\infty$ algebra structure on $V$.
Consider the tangent space to the space of these maps at some point, denote
this tangent space by $\Der(M^\mb,\End(V))$. Then we can construct from the
differential $\partial$ in $M^\mb$ an odd vector field on the space
$\Der(M^\mb,\End(V))[-1]$ as we explain in Section 3. In this case of associative algebras, this construction gives
exactly the Hochschild complex with the Gerstenhaber bracket. This is the
alternative construction of the deformation Lie (in general, $L_\infty$)
algebra.

The initial problem from which this paper was grown, is the deformation
theory (in the classical, commutative sense) of the (co)associative
bialgebras. We tried to construct an $L_\infty$ algebra governing the
deformations of an (co)associative bialgebra. We supposed that the
underlying complex of this $L_\infty$ algebra is quasi-isomorphic to the
Gerstenhaber-Schack complex of the (co)associative bialgebra. We had in mind
the operadic construction with the Stasheff associahedrons described above.
Then we constructed a compactification of some spaces $\K_m^n$ introduced by
Maxim Kontsevich. It turned out that the boundary strata of this, a very natural, compatification,
are the products of not only the spaces $\K_{m^\prime}^{n^\prime}$, but some
spaces with multiindices $\K_{m_1,\dots,m_{\ell_1}}^{n_1,\dots, n_{\ell_2}}$
generalizing the space $\K_m^n$. Then it turned out that the algebraic
operations under these spaces form not a PROP as was expected but a new
algebraic concept called here a CROC. We construct a CROC $\End(V)$ for a
vector space $V$. Then we have a construction which is a direct
generalization of the construction with the minimal model $M^\mb$ in the
Stasheff case, but here this construction gives us an {\it associative dg
algebra}, not a Lie dg algebra. Then we interpret this associative algebra
as the associative algebra governing the {\it non-commutative} deformations
of the initial bialgebra. In the same way we formulate a concept of a
non-commutative homotopical bialgebra.

We have analogs of all these constructions in the well-understood case of
associative algebras. In this case, we also construct an associative algebra
governing the non-commutative deformations. According to the Quillen duality
philosophy, we should have "a map" $\g^\mb\to (\T(A^*[1]),\delta)$ where the
latter is the associative algebra we constructed, and $\g^\mb$ is the
Hochschild complex with the Gerstenhaber bracket. We construct such a map
explicitly.

It happens, however, that in the case of
(co)associative bialgebras we {\it can not} construct such a map. The only
what remains is to use again the Quillen duality. We first consider the
Quillen dual to the constructed associative algebra, then take its quotient
by the commutant. It is a dg commutative algebra, conjecturaly the algebra
of functions on the extended commutative moduli space in the formal
neighborhood of the point corresponding to the initial bialgebra. Then we
construct the Quillen dual dg Lie algebra which conjecturally governs the
deformations of the bialgebra.

\section{{\tt The space $\K(m,n)$ and its compactification}}
The space $\K(m,n)$ we consider here and its compactification play the same role in the deformation theory
of associative bialgebras as the Stasheff polyhedra play in
the deformation theory of associative algebras. Here we define this space and construct its
compactification $\overline{\K(m,n)}$ which is a (compact) manifold with corners.
In the next Sections we define an $L_\infty$-algebra structure on a
complex quasi-isomorphic to the
Gerstenhaber-Schack complex of an associative bialgebra using this
compactification. The spaces $\K(m,n)$ were constructed by Maxim
Kontsevich.

\subsection{{\tt The space $\K(m,n)$.}}
First define the space $\Conf (m,n)$. By definition, $m,n\ge 1$, $m+n\ge 3$,
and
\begin{multline}\label{eq1.1}
\Conf (m,n)=\{p_1,\dots, p_m\in \mathbb{R}^{(1)}, p_i<p_j\ \  for\ \  i<j;\\
q_1,\dots,q_n\in\mathbb{R}^{(2)}, q_i<q_j\ \ for\ \ i<j\}
\end{multline}
Here we denote by $\mathbb{R}^{(1)}$ and by $\mathbb{R}^{(2)}$ two different
copies of a real line $\mathbb{R}$.

Next, define a 3-dimensional group $G^3$ acting on $\Conf(m,n)$. This group
is a semidirect product $G^3=\mathbb{R}^2\ltimes\mathbb{R}_+$ (here
$\mathbb{R}_+=\{x\in\mathbb{R}, x>0\}$) with the following group law:
\begin{equation}\label{eq1.2}
(a,b,\lambda)\circ(a^{\prime},b^{\prime},\lambda ^{\prime})=
(\lambda ^{\prime} a+a^{\prime},(\lambda ^{\prime})^{-1} b+b^{\prime},\lambda\lambda ^{\prime} )
\end{equation}
where $a,b,a^\prime ,b^\prime\in\mathbb{R}, \lambda,\lambda^\prime\in\mathbb{R}_+$.
This group acts on the space $\Conf(m,n)$ as
\begin{multline}\label{eq1.3}
(a,b,\lambda)\cdot (p_1,\dots,p_m;q_1,\dots,q_n)=
(\lambda p_1+a,\dots,\lambda
p_m+a;\lambda^{-1}q_1+b,\dots,\lambda^{-1}q_n+b)
\end{multline}
In other words, we have two independent shifts on $\mathbb{R}^{(1)}$ and
$\mathbb{R}^{(2)}$ (by $a$ and $b$), and $\mathbb{R}_+$ dilatates
$\mathbb{R}^{(1)}$ by $\lambda$ and dilatates $\mathbb{R}^{(2)}$ by $\lambda
^{-1}$.

In our conditions $m,n\ge 1, m+n\ge 3$, the group $G^3$ acts on $\Conf(m,n)$
freely. Denote by $\K(m,n)$ the quotient-space. It is a smooth manifold of
dimension $m+n-3$.
\subsubsection{{\tt Example}}
Let $m=n=2$. Then the space $\K(2,2)$ is 1-dimensional. It is easy to see
that $(p_2-p_1)\cdot (q_2-q_1)$ is preserved by the action of $G^3$, and it
is the only invariant of the $G^3$-action on $\K(2,2)$. Therefore, $\K(2,2)\simeq
\mathbb{R}_+$. There are two "limit" configurations: $(p_2-p_1)\cdot
(q_2-q_1)\rightarrow 0$ and $(p_2-p_1)\cdot (q_2-q_1)\rightarrow \infty$.
Therefore, the compactification $\overline{\K(2,2)}\simeq [0,1]$.

The main trouble in the problem of constructing the (right) compactification
$\overline{\K(m,n)}$ in the general case is that the space $\K(m,n)$ is not
compact, and the points can move away from each other on infinite distances.
Moreover, all these infinities are not the same, in particular, $\infty$ and
$\infty^2$ are different infinities when they occur in the same
configuration. See the following example.

\subsubsection{{\tt Example}}
Consider the space $\K(1,n)$, $n\ge 2$. Consider the following "limit"
configuration in $\K(n,1)$: first two points $q_1,q_2\in \mathbb{R}^{(2)}$
are in a finite distance from each other; the point $q_3$ is in the distance
$\infty$ from $q_1,q_2$; the point $q_4$ is in the distance $\infty^2$ from
$q_3$; $q_5$ is in the distance $\infty^3$ from $q_4$, and so on. We will
see in the next Subsection that the space of all such configurations in
$\K(1,n)$ has dimension 0.

Now we are going to construct the compactification $\overline{\K(m,n)}$ in the
general case.

\subsection{{\tt The compactification $\overline{\K(m,n)}$.}}
\subsubsection{{\tt The space
$\K^{m_1,\dots,m_{\ell_1}}_{n_1,\dots,n_{\ell_2}}$.}}
Define first the space $\K^{m_1,\dots,m_{\ell_1}}_{n_1,\dots,n_{\ell_2}}$ of
dimension $\sum_{i=1}^{\ell_1}m_i+\sum_{i=1}^{\ell_2}n_i-\ell_1-\ell_2-1$
(here $m_i,n_i\ge 1$ and
$\sum_{i=1}^{\ell_1}m_i+\sum_{i=1}^{\ell_2}n_i\ge\ell_1+\ell_2+1$).
The space $\K(m,n)$ a particular case of these spaces, when
$\ell_1=\ell_2=1$, $\K(m,n)=\K^m_n$. Fist define the space
$\Conf^{m_1,\dots,m_{\ell_1}}_{n_1,\dots,n_{\ell_2}}$. By definition,
\begin{multline}\label{eq1.4}
\Conf^{m_1,\dots,m_{\ell_1}}_{n_1,\dots,n_{\ell_2}}=\\
\{p^1_1,\dots,p^1_{m_1}\in\mathbb{R}^{(1,1)},
p^2_1,\dots,p^2_{m_2}\in\mathbb{R}^{(1,2)},\dots,
p^{\ell_1}_1,\dots,p^{\ell_1}_{m_{\ell_1}}\in\mathbb{R}^{(1,\ell_1)};\\
q^1_1,\dots,q^1_{n_1}\in\mathbb{R}^{(2,1)},
q^2_1,\dots,q^2_{n_2}\in\mathbb{R}^{(2,2)}\dots,
q^{\ell_2}_1,\dots,q^{\ell_2}_{n_{\ell_2}}\in\mathbb{R}^{(2,\ell_2)}|\\
p^j_{i_1}<p^j_{i_2}\ \ for\ \ i_1<i_2;q^j_{i_1}<q^j_{i_2}\ \ for\ \
i_1<i_2\}
\end{multline}
Here $\mathbb{R}^{(i,j)}$ are copies of the real line $\mathbb{R}$.
Now we have an $\ell_1+\ell_2+1$-dimensional group $G^{\ell_1,\ell_2,1}$
acting on $\Conf^{m_1,\dots,m_{\ell_1}}_{n_1,\dots,n_{\ell_2}}$.
It contains $\ell_1+\ell_2$ independent shifts
$$
p_i^j\mapsto p_i^j+a_j, i=1,\dots, m_j, a_j\in\mathbb{R};
q_i^j\mapsto q_i^j+b_j, i=1,\dots, n_j, b_j\in\mathbb{R}
$$
and {\it one } dilatation
$$
p_i^j\mapsto \lambda\cdot p_i^j\ \ for\ \ all\ \ i,j;q_i^j\mapsto
\lambda^{-1}\cdot q_i^j\ \ for\ \ all\ \ i,j.
$$
This group is isomorphic to $\mathbb{R}^{\ell_1+\ell_2}\ltimes\mathbb{R}_+$.
We say that the lines
$\mathbb{R}^{(1,1)},\mathbb{R}^{(1,2)},\dots,\mathbb{R}^{(1,\ell_1)}$
(corresponding to the factor $\lambda$) are the lines of the first type, and
the lines $\mathbb{R}^{(2,1)},\mathbb{R}^{(2,2)},\dots,\mathbb{R}^{(2,\ell_2)}$
(corresponding to the factor $\lambda^{-1}$) are the lines of the second
type.

Denote
\begin{equation}\label{eq1.5}
\K^{m_1,\dots,m_{\ell_1}}_{n_1,\dots,n_{\ell_2}}=\Conf
^{m_1,\dots,m_{\ell_1}}_{n_1,\dots,n_{\ell_2}}/G^{\ell_1,\ell_2,1}
\end{equation}

The strata of the compactification $\overline{\K(m,n)}$ constructed below
are direct products of spaces
$\K^{m_1,\dots,m_{\ell_1}}_{n_1,\dots,n_{\ell_2}}$ for different
$\ell_1,\ell_2,m_i,n_j$.
\subsubsection{{\tt The Construction.}}
We define here a stratified manifold $\overline{\K(m,n)}$.

Suppose we have a "limit" configuration $\sigma\in \K(m,n)$; "limit" here
means that some distances $|p_i-p_j|$ or $|q_i-q_j|$ are infinitely small or
infinitely large. In the sequel we say "equal to 0" and "equal to $\infty$"
in these cases.
The non-limit configurations form the maximal open stratum
isomorphic to $\K(m,n)$.

Each limit configuration belongs to a unique stratum of the form
\begin{equation}\label{eq1.6}
\K(\sigma
)\simeq \K^{m_1,\dots,m_{\ell_1}}_{n_1,\dots,n_{\ell_2}}\times
\K^{m^{\prime}_1,\dots,m^{\prime}
_{\ell^{\prime}_1}}_{n^{\prime}_1,\dots,n^{\prime}_{\ell^{\prime}_2}}\times\dots
\end{equation}
which we are going to describe.

Consider the set of all  (finite, infinitely small, or infinitely large)
parameters $\lambda$ such that the image $\sigma_\lambda$ of $\sigma$ after
the application of the element $(0,0,\lambda)\in G^3$obeys the following
property (*):
\begin{quote}
{\tt In $\sigma_\lambda$ {\bf either} there exist at least two points $p_i$ and
$p_{i+1}$ in $\mathbb{R}^{(1)}$ such that the distance $|p_{i+1}-p_i|$ is {\bf finite } (nor
infinitely small neither infinitely large) {\bf or} at least two points
$q_j, q_{j+1}$ in $\mathbb{R}^{(2)}$ with the same property for
$|q_{j+1}-q_j|$.}
\end{quote}

We say that two parameters $\lambda_1 ,\lambda_2$ obeying the property (*)
are equivalent, if the ratio $\frac{\lambda_1}{\lambda_2}$ is finite and not
infinitesimally small.

Denote by $S(\sigma)$ the set of the equivalence classes of the parameters
$\lambda$ obeying the property (*) for the configuration $\sigma$. The set
$S(\sigma)$clearly is not empty for any $\sigma$, and the condition $\sharp
S(\sigma)=1$ is equivalent that $\sigma$ is a non-limit configuration. For a
limit $\sigma$, $\sharp S(\sigma)>1$. It is clear that for any $\sigma$ the
set $S(\sigma)$ is finite.
\begin{example}
Consider the configuration $\sigma$ from Example~1.1.2. We have the space
$\K(1,n)$, and $|q_2-q_1|$ is finite, $|q_3-q_2|\sim\infty$, $|q_4-q_3|\sim
\infty^2$, $\dots , |q_n-q_{n-1}|\sim\infty^{n-2}$. We have: $\sharp
S(\sigma)=n-1$. Indeed, roughly speaking,
$\lambda_1=1,\lambda_2=\infty,\dots, \lambda_{n-1}=\infty^{n-2}$obey the
property (*), it is clear that each $\lambda$ obeying the property (*) for
$\sigma$ is equivalent to some $\lambda_i$ from the list above.
\end{example}

\begin{remark}
We should specify what is meant by a limit configuration. Each point on each
line, when it moves, becomes a real-valued function on a real parameter $t$.
Thus, we have functions $p_1(t),\dots, p_m(t);q_1(t),\dots, q_n(t)$. {\it We
suppose that all these functions are Lourent power series in the parameter
$t$}. Then, a limit configuration is this configuration when $t\to\infty$.
It is important that we do not consider some more other functions in $t$
except polynomials in $t,t^{-1}$. In a sense, thus we obtain the minimal
compactification.
\end{remark}
Consider the configuration $(0,0,\lambda)\cdot\sigma=\sigma_\lambda$.
We identify in $\sigma_\lambda$ any two points which are infinitely close to
each other. Then the equivalence classes (under this identification) of the
points on $\mathbb{R}^{(1)}$ in $\sigma_\lambda$ can be uniquely divided to
$\ell_1(\lambda)$ groups
$$
\{p_1^1,\dots,p_{m_1(\lambda)}^1\},\{p_1^2,\dots,p_{m_2(\lambda)}^2\},\dots,
\{p_1^{\ell_1(\lambda)},\dots,p_{m_{\ell_1(\lambda)}(\lambda)}^{\ell_1(\lambda)}\}
$$
of points standing in turn such that inside each group all the distances
between the points are finite (and nonzero, because we have collapsed the
points infinitely close to each other), and the distances between different
groups are $\infty$. Analogously, we divide the points on $\mathbb{R}^{(2)}$
in $\sigma_\lambda$ to the $\ell_2(\lambda)$ groups
$$
\{q_1^1,\dots,p_{n_1(\lambda)}^1\},\{q_1^2,\dots,q_{n_2(\lambda)}^2\},\dots,
\{q_1^{\ell_2(\lambda)},\dots,q_{n_{\ell_2(\lambda)}(\lambda)}^{\ell_2(\lambda)}\}
$$
by the same way.

We associate with the element $\lambda\in S(\sigma)$ the space
\begin{equation}\label{eq1.7}
\K(\sigma_\lambda)=\K^{m_1(\lambda),\dots,m_{\ell_1(\lambda)
}(\lambda)}_{n_1(\lambda),\dots,n_{\ell_2(\lambda)}(\lambda)}
\end{equation}

\begin{definition}
We say that two limit configurations $\sigma_1, \sigma_2$ are equivalent
iff:
\begin{itemize}
\item[(i)] the sets $S(\sigma_1)$ and $S(\sigma_2)$ are {\bf coincide}
(it means that $\sharp S(\sigma_1)=\sharp S(\sigma_2)$, and we can choose
representatives $\lambda_1,\dots,\lambda_{\sharp S(\sigma_1)}$ for
$S(\sigma_1)$, and representatives
$\lambda^{\prime}_1,\dots,\lambda^{\prime}_{\sharp S(\sigma_2)}$ for
$S(\sigma_2)$ such that there exists $\lambda_{tot}$ (finite, infinitely
small, or infinitely large) such that $\lambda^{\prime}_i=\lambda_{tot}
\cdot\lambda_i$ for each $i$);
\item[(ii)] the sets
$\{\{\ell_1(\lambda)\},\{\ell_2(\lambda)\},\{m_1(\lambda),\dots,m_{\ell_1(\lambda)}(\lambda)\},
\{n_1(\lambda),\dots,n_{\ell_2(\lambda)}(\lambda)\}|\lambda\in S(\sigma_i)\}$
coincide as ordered sets under the identification with $\lambda_{tot}$
described in (i).
\end{itemize}
\end{definition}

In other words, two limit configurations are equivalent if their
"topological types" coincide after the application of some
$(0,0,\lambda_{tot})\in G^3$ for finite, infinitely small, or infinitely
large $\lambda_{tot}$.

\begin{theorem*}
\begin{itemize}
\item[(i)] For a configuration $\sigma$ (limit or non-limit) with $m$ points
on $\mathbb{R}^{(1)}$ and $n$ points on $\mathbb{R}^{(2)}$, the space of all
configurations equivalent to $\sigma$ is homeomorphic to
\begin{equation}\label{eq1.10}
\K(\sigma)=\prod_{[\lambda]\in S(\sigma)}\K(\sigma_\lambda)
\end{equation}
(see (\ref{eq1.7}) for the definition of $\K(\sigma_\lambda)$);
\item[(ii)]
\begin{equation}\label{eq1.11}
 \overline{\K(m,n)}:=\bigsqcup_{\begin{subarray}{c}\text{all equiv.} \\ \text{classes
of}\\  \text{config.\ $\sigma$} \end{subarray}}\K(\sigma)
\end{equation}
defines a compactification of the space $\K(m,n)$; it is a manifold with
corners.
\end{itemize}
\begin{proof}
it is clear. The only what we want to notice is that the distances between
the groups $\{p_1^1,\dots,p_{m_1(\lambda)}^1\},\dots$are $\infty$,therefore,
their positions on the lines are defined up to shifts, for a shift for each
group. It motivates our definition of the space $\K(\sigma_\lambda)$.
\end{proof}
\end{theorem*}
\begin{corollary}
For each $\sigma$, $\dim \K(\sigma)\le m+n-3$.
\qed
\end{corollary}

Although the Corollary follows from Theorem above, it is instructive to give
here a straightforward proof of it. We want to prove that $\dim
\K(\sigma)\le \dim\K(m,n)=m+n-3$, and the equality holds only for a
non-limit configuration $\sigma$.

The set $S(\sigma)$ is naturally ordered by the numbers
$\lambda\in\mathbb{R}_+$ representing the equivalence class. Consider the
minimal element $[\lambda]\in S(\sigma)$. Then there are no 0 distances
between points on $\mathbb{R}^{(2)}$ in $\sigma_\lambda$. Then the points in
$\mathbb{R}^{(2)}$ in the configuration $\sigma_\lambda$ are divided to the
groups $\{q_1^1,\dots,q_{n_1(\lambda)}^1\},\dots,
\{q_1^{\ell_2(\lambda)},\dots,
q_{n_{\ell_2(\lambda)}(\lambda)}^{\ell_2(\lambda)}\}$, and also the points
on $\mathbb{R}^{(1)}$ are divided to the groups (but there can occur 0
distances). We will consider only that part of $\dim\K(\sigma)$ which is
contributed by the points on $\mathbb{R}^{(2)}$, the contribution of points
on $\mathbb{R}^{(1)}$ is analogous. We denote this dimension by $\dim_2$.
Thus, for the minimal $\lambda\in S(\sigma)$,
$\dim_2\K(\sigma_\lambda)=n_1(\lambda)+\dots+n_{\ell_2(\lambda)}(\lambda)-\ell_2(\lambda)$.
Also we have $\dim
\K(\sigma_\lambda)=\dim_1\K(\sigma_\lambda)+\dim_2\K(\sigma_\lambda)-1$ (the
last $-1$ because of the action of $(0,0,\mathbb{R}_+)\in
G^{\ell_1,\ell_2,1}$). Now consider the next (in the sense of the canonical
ordering) element $\lambda^{\prime}\in S(\sigma)$. It is clear that
$\frac{\lambda^{\prime}}{\lambda}=\infty$. Each group of points
$\{q_1^i,\dots,q_{n_i(\lambda)}^i\}$ will be collapsed to a point
$\overline{q_{n_i}}$ in $\sigma_{\lambda^{\prime}}$, and only these groups will
be collapsed because of our choice of $\lambda^{\prime}$. Then, in
$\sigma_{\lambda^\prime}$ the points $\overline{q_1},\dots,
\overline{q_{\ell_2(\lambda)}}$ are divided to $\ell_2(\lambda^\prime)$ groups,
and
$$
\dim_2\K(\sigma_{\lambda^\prime})=\ell_2(\lambda)-\ell_2(\lambda^\prime)
$$
We see that
$$
\dim_2\K(\sigma_\lambda)+\dim_2\K(\sigma_{\lambda^\prime})=
n_1(\lambda)+\dots+n_{\ell_2(\lambda)}-\ell_2(\lambda^\prime)=n-\ell_2(\lambda^\prime)
$$
We can show analogously, that after $d$ steps,
$$
\dim_2\K(\sigma_{\lambda_1})\times\dots\times\K(\sigma_{\lambda_d})=
n-\ell_2(\lambda_d)
$$
(here $\lambda_1<\dots<\lambda_d$ are the first lowest $d$ elements in
$S(\sigma)$).
For the maximal $\lambda_{max}\in S(\sigma)$, $\ell_2(\lambda_{max})=1$
(there is the only 1 group of points). Therefore, $\dim_2\K(\sigma)=n-1$.
It proves the following
\begin{proposition}
$\dim\K(\sigma)=m+n-\sharp S(\sigma)-2$
\qed
\end{proposition}
\begin{remark}
In the case $m=1$ the space $\overline{\K(1,n)}$ is {\it NOT} the Stasheff polyhedron.
We can see it immediately. In particular, in our compactification it is
important the "relative velocity" with which points move close to each
other. In fact, in the Stasheff compactification it is irrelevant. This fact
could hint us that the algebraic structures behind our compactification is
not anymore some usual structures like OPERADs and PROPs.
\end{remark}
\subsubsection{{\tt Examples}}
\paragraph{{\tt Example}}
We already know that the space $\K(2,2)$ has 2 different limit
configurations: they are $\sigma_1$, when $|p_2-p_1|\cdot|q_2-q_1|\sim 0$,
and $\sigma_2$, when $|p_2-p_1|\cdot|q_2-q_1|\sim \infty$ (see Example
1.1.1).We have:
\begin{equation}\label{eq1.14}
\K(\sigma_1)=\K_2^1\times\K_1^2
\end{equation}
\begin{equation}\label{eq1.15}
\K(\sigma_2)=\K^{1,1}_2\times\K^2_{1,1}
\end{equation}
In both cases $\dim\K(\sigma_i)=0$.
\paragraph{{\tt Example}}
Consider the limit configuration $\sigma$ in the space $\overline{\K(1,n)}$
described in the Example 1.1.2. Recall that for this configuration $\sigma$
we have: $|q_2-q_1|$ is finite, $|q_3-q_2|\sim\infty$,
$|q_4-q_3|\sim\infty^2$,..., $|q_n-q_{n-1}|\sim\infty^{n-2}$.
We have: $\sharp S(\sigma)=n-1$, and therefore $\dim\K(\sigma)=0$ by
Proposition 1.2.2. We have:
$$
\K(\sigma)=\K^1_{2,1,\dots,1 \\(n-2\ of\ 1's)}\times\K^1_{2,1,\dots,1
(n-3\ of\ 1's)}\times\dots \times\K_2^1.
$$
\subsection{{\tt Strata of codimension 1}}
Here we describe all strata of codimension 1 in $\overline{\K(m,n)}$. These
strata are very important in the next Sections where we construct an
$L-\infty$-algebra structure on the deformation complex of an associative
bialgebra.

A typical stratum of codimension 1 is drawn in Figure~1.
\sevafigc{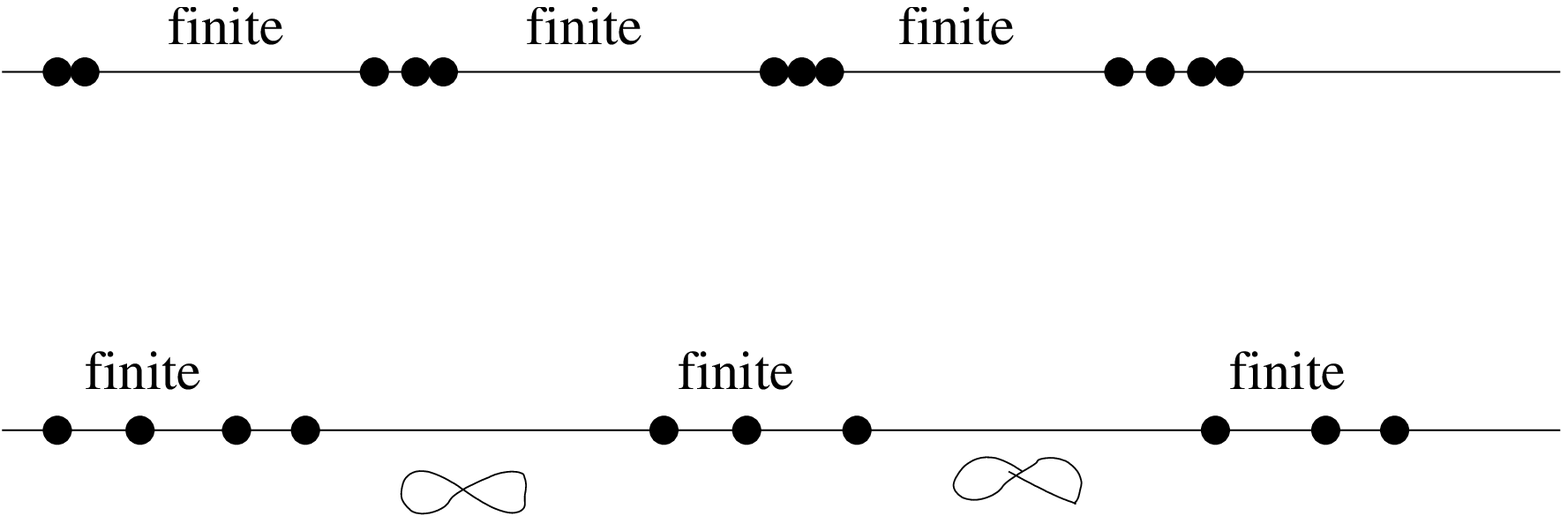}{90mm}{0}{A typical stratum of codimension 1}

Consider numbers $\ell_1$, $\ell_2$, $\{m_1,\dots,m_{\ell_1}\}$, $\{n_1,\dots,
n_{\ell_2}\}$ satisfying $\sum_{i=1}^{\ell_1}m_i=m,\
\sum_{i=1}^{\ell_2}n_i=n$.
Consider the limit configuration $\sigma\in\overline{\K(m,n)}$
in which the points on $\mathbb{R}^{(1)}$ are divided to $\ell_1$ groups,
$m_i$ points in the $i$-th group, and the points inside each group are
infinitely close to each other (in order $\varepsilon$); the points on
$\mathbb{R}^{(2)}$ are divided to $\ell_2$ groups, $n_j$ points in the
$j$-th group, and points inside each group are in finite distance from each
other, and the points of different groups are infinitely far from each other
(in order $\frac1{\varepsilon}$). It is clear that $\sharp S(\sigma)=2$ (representatives in $S(\sigma)$ are
$\lambda=1$ and $\lambda=\frac1{\varepsilon}$), and
it is the most general configuration with $\sharp S(\sigma)=2$. Then it
follows from Proposition 1.2.2 that these configurations $\sigma$ exhaust
all strata of codimension 1.

\subsection{\tt{The space $\overline{\K_{m_1,\dots , m_{\ell_1}}^{n_1 , \dots,n_{\ell_2}}}$}}

Here we construct the compactification
$\overline{\K_{m_1,\dots,m_{\ell_1}}^{n_1,\dots ,n_{\ell_2}}}$ of the space
$\K_{m_1,\dots,m_{\ell_1}}^{n_1,\dots,n_{\ell_2}}$ which we will need in
the sequel.

The construction is analogous to the construction of $\overline{\K(m,n)}$
above. We associate with a configuration in ${\K_{m_1,\dots , m_{\ell_1}}^{n_1 , \dots,n_{\ell_2}}}$
a limit configuration in $\K(\sum_{i=1}^{\ell_1} m_i,\sum_{j=1}^{\ell_2}n_j)$. Namely, we divide
the $\sum_{i=1}^{\ell_1}m_i$ points into $\ell_1$ groups with $m_i$ points
in the $i$-th group, and do the same with the second line. We suppose that
the distances between groups are $\infty^N$ where $N\gg 0$. In other words,
we suppose that the distances between the groups are infinitely large
comparably with all other infinities in the (sub)limit configurations we
consider. Then the previous construction can be easily generalized to this
case.

\section{{\tt The concept of CROC}}

Here we introduce the concept of {\it CROC}. This notion formalizes the
admissible operations on an algebraic structure, like OPERADs and PROPs.
We tried to formalize operations which we have among the spaces
$\overline{\K_{m_1,\dots,m_{\ell_1}}^{n_1,\dots,n_{\ell_2}}}$ which is an
example of a topological CROC. It turns out that we can describe the
associative bialgebras as algebras over some CROC, and this CROC of
associative bialgebras has a very natural simplicial free minimal model,
analogous to the Stasheff construction in the case of associative algebras.
This simplicial free resolution is formed from the chain complexes of the
spaces $\overline{\K_{m_1,\dots,m_{\ell_1}}^{n_1,\dots,n_{\ell_2}}}$.

\subsection{{\tt The Definition}}
\begin{definition}
{\tt A preCROC of vector spaces is a collection of vector spaces
$F_{m_1,\dots,m_{\ell_1}}^{n_1,\dots,n_{\ell_2}}$, $m_i,n_j\ge 1$ with a
left action of the product of symmetric groups
$\Sigma_{m_1}\times\dots\times\Sigma_{m_{\ell_1}}$ and a right action of
$\Sigma_{n_1}\times\dots\times\Sigma_{n_{\ell_2}}$
and with a composition law. In the simplest case, this composition law is a
map
\begin{equation}\label{eqc1.1}
\phi_{m_1,\dots, m_{\ell_1}|\ell_2}^{\ell_1 |n_1,\dots, n_{\ell_2}}\colon F_{m_1,\dots,m_{\ell_1}}^{\ell_2}\otimes
F_{\ell_1}^{n_1,\dots,n_{\ell_2}}\to
F_{m_1+m_2+\dots+m_{\ell_1}}^{n_1+n_2+\dots +n_{\ell_2}}
\end{equation}
In general, we have the composition law
\begin{multline}\label{eqc1.2}
\phi_{m^1_1,\dots, m^1_{\ell_1^1};\dots;m^a_1,\dots,m^a_{\ell_1^a}|
\ell_2^1,\dots,\ell_2^b}^{\ell_1^1,\ell_1^2,\dots,\ell_1^a|n^1_1,\dots,n^1_{\ell_2^1};\dots;
n_1^b,\dots,n^b_{\ell_2^b}}\colon F_{m^1_1,\dots,
m^1_{\ell_1^1},\dots,m^a_1,\dots,m^a_{\ell_1^a}}^{\ell_2^1,\dots,\ell_2^b}\otimes
F^{n^1_1,\dots,n^1_{\ell_2^1},\dots,
n_1^b,\dots,n^b_{\ell_2^b}}_{\ell_1^1,\dots,\ell_1^a}\to\\
F_{\sum_{i=1}^{\ell_1^1}m_i^1,\dots,\sum_{i=1}^{\ell_1^a}m^a_i}^
{\sum_{j=1}^{\ell_2^1}n_j^1,\dots,\sum_{j=1}^{\ell^b_2}n_j^b}
\end{multline}
There are three axioms on this data:
\begin{itemize}
\item[(i)]$F_1^1$ is the 1-dimensional trivial representation of
$\Sigma_1\times\Sigma_1$;
\item[(ii)] the compositions are compatible with the action of symmetric
groups;
\item[(iii)] the natural associativity of the compositions.
\end{itemize}}
\end{definition}

\begin{definition}
{\tt A CROC is a preCROC with the following extra conditions:
\begin{itemize}
\item[(i)]there are the {\it restriction maps} $r_i\colon
F_{m_1,\dots, m_{\ell_1}}^{n_1,\dots, n_{\ell_2}}\to
F_{m_1,\dots,\hat{m_i},\dots, m_{\ell_1}}^{n_1,\dots, n_{\ell_2}}$
and $r^j\colon F_{m_1,\dots, m_{\ell_1}}^{n_1,\dots, n_{\ell_2}}\to
F_{m_1,\dots, m_{\ell_1}}^{n_1,\dots,\hat{n_j},\dots, n_{\ell_2}}$,
\item[(ii)]the restrictions satisfy the natural commutativity of the
compositions.
\end{itemize}}
\end{definition}
We will consider algebras over (pre)CROCs. For this we define the preCROC $\End(V)$
and define an algebra over a preCROC $F$ structure on a vector space $V$ as a
map of preCROCs $F\to\End(V)$. Note that the definition of the preCROC $\End(V)$
is not very straightforward, we do it in the next Subsection.

\subsection{{\tt The preCROC $\End(V)$}}

Let $V$ be a vector space. Here we define the CROC $\End(V)$.
By definition,
\begin{equation}\label{eqc1.3}
\End(V)_{m_1,\dots,m_{\ell_1}}^{n_1,\dots, n_{\ell_2}}=
\bigotimes_{i=1...\ell_1, j=1...\ell_2}\Hom(V^{\otimes m_i},V^{\otimes n_j})
\end{equation}
We should define now the composition and the maps $\vartheta_i$,
$\vartheta^i$. First define the composition.

We define first the simplest composition
$$
\phi_{m_1,\dots, m_{\ell_1}|\ell_2}^{\ell_1 |n_1,\dots, n_{\ell_2}}\colon \End(V)_{m_1,\dots,m_{\ell_1}}^{\ell_2}\otimes
\End(V)_{\ell_1}^{n_1,\dots,n_{\ell_2}}\to
\End(V)_{m_1+m_2+\dots+m_{\ell_1}}^{n_1+n_2+\dots +n_{\ell_2}}
$$
Suppose we have
$$
\Psi_1\in\Hom(V^{\otimes \ell_1},V^{\otimes n_1}),\Psi_2\in\Hom(V^{\otimes \ell_1},V^{\otimes
n_2}),\dots,\Psi_{\ell_2}\in\Hom(V^{\otimes \ell_1},V^{\otimes
n_{\ell_2}})
$$
and
$$
\Theta_1\in\Hom(V^{\otimes
m_1},V^{\otimes \ell_2}),\Theta_2\in\Hom(V^{\otimes
m_2},V^{\otimes \ell_2}),\dots,\Theta_{\ell_1}\in\Hom(V^{\otimes
m_{\ell_1}},V^{\otimes \ell_2}),
$$
we are going to define their composition which belongs to
$\End(V)_{m_1+\dots+m_{\ell_1}}^{n_1+\dots+n_{\ell_2}}$. Denote
$m=m_1+\dots+m_{\ell_1}$, $n=n_1+\dots+m_{\ell_2}$. The construction is
as follows:

First define
\begin{multline}\label{eqc2.6}
F(v_1\otimes\dots\otimes v_m)\colon =
\Theta_1(v_1\otimes\dots\otimes
v_{m_1})\bigotimes\Theta_2(v_{m_1+1}\otimes\dots\otimes
v_{m_1+m_2})\bigotimes\dots\\ \bigotimes
\Psi_{\ell_2}(v_{m_1+\dots+m_{\ell_1-1}+1}\otimes\dots\otimes v_{{m_1}+\dots+{m_{\ell_1}}})
\in V^{\otimes \ell_1\ell_2}
\end{multline}
Now we apply $\{\Theta_i\}$ to this element in $V^{\otimes \ell_1\ell_2}$:
we define an element $G\colon V^{\otimes \ell_1\ell_2}\to V^{\otimes m}$ as follows:
\begin{multline}\label{eqc2.7}
G(v_1\otimes v_2\otimes\dots\otimes v_{\ell_1\ell_2}):=
\Theta_1(v_1\otimes v_{\ell_2+1}\otimes\dots\otimes
v_{\ell_2(\ell_1-1)+1})\bigotimes\\
\Theta_2(v_2\otimes v_{\ell_2+2}\otimes\dots\otimes
v_{\ell_2(\ell_1-1)+2})\bigotimes\dots
\bigotimes\Theta_{\ell_1}(v_{\ell_2}\otimes v_{2\ell_2}\otimes\dots\otimes
v_{\ell_2\ell_1})\in V^{\otimes n}.
\end{multline}
Define now
\begin{equation}\label{eqc2.8}
Q(v_1\otimes \dots\otimes v_m):=G\circ
F(v_1\otimes\dots\otimes v_m)
\end{equation}

By definition, the element $Q$ is the composition
$\phi_{m_1,\dots, m_{\ell_1}|\ell_2}^{\ell_1 |n_1,\dots,
n_{\ell_2}}(\Theta_1\otimes\dots\otimes\Theta_{\ell_1}\otimes\Psi_1\otimes\dots\otimes\Psi_{\ell_2})$

This construction can be easily generalized to the higher compositions
$\phi_{m^1_1,\dots, m^1_{\ell_1^1};\dots;m^a_1,\dots,m^a_{\ell_1^a}|
\ell_2^1,\dots,\ell_2^b}^{\ell_1^1,\ell_1^2,\dots,\ell_1^a|n^1_1,\dots,n^1_{\ell_2^1};\dots;
n_1^b,\dots,n^b_{\ell_2^b}}$ such that the associativity holds. The formulas
are very huge but the pictures behind them are very simple.

We can easily prove that the structure defined in this way is a preCROC.

\begin{definition}
\begin{itemize}
\item[(i)] {\tt Let $F$ be a preCROC. Then an $F$-algebra structure on a vector
space $V$ is a map of preCROCs $F\to\End(V)$,
\item[(ii)] Let $F$ be a CROC. Then an $F$-algebra on a vector space $V$ is
a map of preCROCs $\Upsilon\colon F\to\End(V)$ such that
$\Upsilon(F_1^1)=\Id\in\Hom(V,V)$ and which is compatible with the
restrictions on $F$ and with the tensor product on $\End(V)$:
\begin{equation}\label{eqcc1}
\Upsilon(F_{m_1,\dots, m_{\ell_1}}^{n_1,\dots, n_{\ell_1}})=\Upsilon(r_i(F_{m_1,\dots, m_{\ell_1}}^{n_1,\dots,
n_{\ell_1}}))\otimes\Upsilon(\bar{r_i}(F_{m_1,\dots, m_{\ell_1}}^{n_1,\dots,
n_{\ell_1}}))
\end{equation}
where $\bar{r_i}\colon F_{m_1,\dots, m_{\ell_1}}^{n_1,\dots, n_{\ell_1}}\to
F_{m_i}^{n_1,\dots, n_{\ell_1}}$ is the "complementary" restriction which is
defined as the composition of the corresponding restrictions. The same should be true for
$r^j$.}
\end{itemize}
\end{definition}
\subsection{{\tt The CROC of associative bialgebras}}
In this Subsection we define the CROC $\mathrm{Assoc}$ of associative
bialgebras and show that a map of CROCs $\mathrm{Assoc}\to\End(V)$ is the
same that an associative bialgebra structure on $V$ (here by an associative
bialgebra we mean a Hopf algebra without the unit, the counit, and the
antipode, that is, it has an associative product, a coassociative coproduct,
which are compatible).

By definition, $\mathrm{Assoc}_{m_1,\dots, m_{\ell_1}}^{n_1,\dots,
n_{\ell_2}}=\prod_{i=1...\ell_1,j=1...\ell_2}\Sigma_{m_i}\times\Sigma_{n_j}$
where $\Sigma_i$ is the symmetric group. Now we are going to define the compositions.

For this, for any associative bialgebra $V$ we associate with an element in
$\Sigma_{m_i}\times\Sigma_{n_j}$ the following element
in $\Hom(V^{\otimes m_i},V^{\otimes n_j})$:
$\Psi_{i,j}\colon v_1\otimes\dots\otimes v_{m_i}\mapsto
\sigma_{i,j}^{(2)}\circ
\Delta^{n_j-1}(v_{\sigma_{i,j}^{(1)}1}\star\dots\star
v_{\sigma_{i,j}^{(1)}m_i})$
where we denote by $\sigma_{i,j}^{(1)},\sigma_{i,j}^{(2)}$ the corresponding
permutations from the symmetric groups.
Thus, we attached to each element in $\Sigma_{m_i}\times\Sigma_{n_j}$ an
element in $\Hom(V^{\otimes m_i},V^{\otimes n_j})$ for any bialgebra $V$.
We claim that there exists a unique CROC structure on $\mathrm{Assoc}$ such
that for any bialgebra V the constructed map is a map of CROCs
$\mathrm{Assoc}\to\End(V)$. Indeed, the composition of the corresponding
$\Hom$'s is again a homomorphism of this form because of the associativity,
the coassociativity, and the compatibility of the product with the
coproduct. We can write down the corresponding permutation by an explicit
formula. We do not do that because this formula will not tell us anything
new.

Now we are going to prove the following result:
\begin{theorem*}
Any map of CROCs $\Upsilon\colon\mathrm{Assoc}\to\End(V)$ is equivalent to the map described above for some
bialgebra
structure on $V$.
\begin{proof}
We already shown that any bialgebra structure on $V$, by the definition of
the CROC $\mathrm{Assoc}$, gives a map of CROCs $\mathrm{Assoc}\to\End(V)$.
To prove the reverse statement, first denote by
$\Psi=\Upsilon(\mathrm{Assoc}_2^1)\in\Hom(V^{\otimes 2},V)$ and
$\Delta=\Upsilon(\mathrm{Assoc}_1^2)\in\Hom(V,V^{\otimes 2})$. We want to prove the
associativity for $\Psi$, the coassociativity for $\Delta$, and their
compatibility. To prove say the associativity, consider the maps (CROC's
compositions)
$i_1\colon\mathrm{Assoc}_{1,2}^1\times
\mathrm{Assoc}_2^1\to\mathrm{Assoc}_3^1$ and
$i_2\colon\mathrm{Assoc}_{2,1}^1\times\mathrm{Assoc}_2^1\to\mathrm{Assoc_3^1}$
It is clear that they coincide when applied to the identity elements of the symmetric groups.
Then, as we have a map
of CROCs, the CROC compositions of their images also should coincide. We
prove the coassociativity in the same way. To prove the compatibility, note
that the two maps
$t_1\colon\mathrm{Assoc}_2^1\times\mathrm{Assoc}_1^2\to\mathrm{Assoc}_2^2$ and
$t_2\colon\mathrm{Assoc}_{1,1}^2\times\mathrm{Assoc}_2^{1,1}\to\mathrm{Assoc}_2^2$
coincide on the identity elements of the symmetric groups. Then, using the
factorization (the property (ii) of the Definition above) we get the claim.
\end{proof}
\end{theorem*}

\subsection{{\tt A free resolution of the CROC $\mathrm{Assoc}$}}
Consider the direct sum of all chain complexes
$\aleph=\bigoplus_{m_1,...,m_{\ell_1};,n_1,...n_{\ell_2}}C_\mb\overline{
K_{m_1,\dots,m_{\ell_1}}^{n_1,\dots,n_{\ell_2}}}$. As $\{\overline{
K_{m_1,\dots,m_{\ell_1}}^{n_1,\dots,n_{\ell_2}}}\}$ form a topological
CROC, $\aleph$ is a CROC of graded vector spaces, namely,
$\aleph_{m_1,\dots,m_{\ell_1}}^{n_1,\dots,n_{\ell_2}}=C_\mb\overline{
K_{m_1,\dots,m_{\ell_1}}^{n_1,\dots,n_{\ell_2}}}$. Moreover, we have a
differential (the chain differential) on this dg CROC. We can prove the
following theorem:
\begin{theorem*}
\begin{itemize}
\item[(i)] The CROC $\aleph$, as a dg CROC, is free,
\item[(ii)] The cohomology of the CROC $\aleph$ is isomorphic to the CROC
$\mathrm{Assoc}$.
\end{itemize}
\begin{proof}
(i) is clear, (ii) follows from the fact that all spaces $\overline{
K_{m_1,\dots,m_{\ell_1}}^{n_1,\dots,n_{\ell_2}}}$ are contractible.
\end{proof}
\end{theorem*}

\subsection{{\tt The algebra $\mho$}}
Here we define our main object--an associative algebra $\mho$.
As a vector space,
\begin{equation}\label{eqcr1}
\mho=\bigoplus_{m_1,...,m_{\ell_1}\ge 1, n_1,...,n_{\ell_2}\ge 1}
\left(\bigotimes_{1\le i\le \ell_1, 1\le j\le\ell_2}\Hom(V^{\otimes
m_i},V^{\otimes n_j})\right)[-\sum m_i-\sum n_j+\ell_1+\ell_2]
\end{equation}
Note that the grading is compatible with the grading in the
Gerstenhaber-Schack complex. Now we define an associative product on $\mho$.
Let $\Psi_1\in\left(\bigotimes_{1\le i\le \ell_1, 1\le j\le\ell_2}\Hom(V^{\otimes
m_i},V^{\otimes n_j})\right)[-\sum m_i-\sum n_j+\ell_1+\ell_2]$ and
$\Psi_2\in\left(\bigotimes_{1\le i\le \ell^\prime_1, 1\le j\le\ell^\prime_2}\Hom(V^{\otimes
m_i^\prime},V^{\otimes n^\prime_j})\right)[-\sum m^\prime_i-\sum
n^\prime_j+\ell^\prime_1+\ell^\prime_2]$. Define their product as the
composition in the preCROC $\End(V)$ (it is 0, if the corresponding composition
$\phi_{\dots}^{\dots}$ in the preCROC $\End(V)$ is 0) {\it up to a sign}.
This sign is defined geometrically from the boundary operator in $\aleph$ as
follows.

For $\Psi_1$ and $\Psi_2$ as above, if their product is nonzero, there exist
{\it a unique} stratum of codimension 1 in {\it a unique} space
$\overline{K_{\dots}^{\dots}}$ which is up to {\it a sign} the space
$\K_{m_1,\dots,m_{\ell_1}}^{n_1,\dots,n_{\ell_2}}\times
\K_{m_1^\prime,\dots, m^\prime_{\ell^\prime_1}}^{n_1^\prime,\dots,
n^\prime_{\ell_2^\prime}}$.

By definition, this sign is equal to the sign in the $\Psi_1\circ \Psi_2$
before their product in the preCROC $\End(V)$.

\begin{theorem*}
The product in $\mho$, defined in this way, is associative.
\begin{proof}
It is clear that the product is associative up to a sign, because the
product in the preCROC $\End(V)$ is associative. Only what we need to check
are the signs.

We have the canonical projection of CROCs $p\colon\aleph\to\mathrm{Assoc}$.
Then, any associative bialgebra is an algebra over the CROC $\aleph$.
Consider the tangent space to the space of  maps of {\it pre}CROCs $\Der(\aleph,\End(V))$ at the point,
corresponding to the bialgebra above. We define a differential and a product
on $\Der(\aleph,\End(V))[-1]$ as follows. Let $\overline{\aleph}$ be the
space of the generators of the free CROC $\aleph$, namely,
$\overline{\aleph}$ consists from the all cells of codimension 0. Then any
element $D\in\Der(\aleph, \End(V)$ is uniquely defined by its restriction
to $\overline{\aleph}$. If we have two derivations $D_1, D_2\in
\Der(\aleph,\End(V))$ we can take the composition
\begin{equation}\label{eqcr2}
\overline{\aleph}\to \overline{\aleph}^{\otimes 2}\to \End(V)^{\otimes 2}\to
\End(V)
\end{equation}
where the first arrow is the chain differential $\partial$, the second is
$D_1\otimes D_2$, and the third is the composition in the CROC $\End(V)$.
It is clear that this definition of the product in $\mho$ coincides with the
definition given above. The advantage of the definition (\ref{eqcr2}) is
that here the signs are specified. But we need to prove that this formula
gives indeed an associative product.

We do it using the equation $\partial^2=0$. Namely, suppose that
$D_1\in\End(V)_{m_1^1,\dots,
m_{\ell_1^1}^1;m_1^2,\dots,m_{\ell_1^2}^2;\dots;m^k_1,\dots,m^M_{\ell^M_1}}^{N}$,
$D_2\in\End(V)_{\ell_1^1,\dots,\ell_1^M}^{\ell_2^1,\dots,\ell_2^N}$,
$D_3\in
\End(V)_M^{n_1^1,\dots,n^1_{\ell_2^1};n^2_1,\dots,n^2_{\ell_2^2};\dots;
n^N_1,\dots,n^N_{\ell_2^N}}$. We want to write down explicitly what follows
from the equation
\begin{equation}\label{eqcr3}
\partial^2(\K_{m_1^1+\dots+m^1_{\ell_1^1}+\dots+m^M_1+\dots+m^M_{\ell^M_1}}^{n^1_1+
n^1_2+\dots+n^1_{\ell_2^1}+\dots+n_1^N+\dots+n_{\ell_2^N}^N})(D_1\otimes D_2\otimes D_3)=0
\end{equation}
One can show that the only interesting boundaries (of codimension 1) in $\partial(\K_{m_1^1+\dots+m^1_{\ell_1^1}+\dots+m^M_1+\dots+m^M_{\ell^M_1}}^{n^1_1+
n^1_2+\dots+n^1_{\ell_2^1}+\dots+n_1^N+\dots+n_{\ell_2^N}^N})$ are:
\begin{equation}\label{eqcr4}
\partial_1=\pm\K_{m_1^1+\dots+m_{\ell_1^1}^1,m_1^2+\dots+m_{\ell_1^2}^2,\dots,m^M_1+
\dots+m^M_{\ell^M_1}}^{\ell_2^1+\dots+\ell_2^N}
\times\K^{n_1^1,\dots,n^1_{\ell_2^1};\dots;n^N_1+\dots+n^N_{\ell^N_2}}_M
\end{equation}
and
\begin{equation}\label{eqcr5}
\partial_2=\pm\K^N_{m^1_1,\dots,m^1_{\ell_1^1};\dots;m^M_1,\dots,m^M_{\ell_1^M}}\times
\K_{\ell_1^1,\dots,\ell_1^M}^{n_1^1+\dots n^1_{\ell_2^1},\dots,n_1^N+\dots
n^N_{\ell_2^N}}
\end{equation}
The boundary of the first factor in (\ref{eqcr4}) contains the term $\pm
\K^N_{m_1^1,\dots,m_{\ell_1^1}^1;\dots;m^M_1,\dots, m^M_{\ell^M_1}}\times
\K_{\ell_1^1,\dots,\ell_1^M}^{\ell_1^1,\dots,\ell_2^N}$ and the
second factor in (\ref{eqcr5}) contains the term
$\pm\K_{\ell_1^1,\dots,\ell_1^M}^{\ell_2^1,\dots,\ell_1^N}\times
\K_M^{n_1^1,\dots, n_{\ell_2^1}^1;\dots;n^N_1,\dots,n^N_{\ell^N_2}}$.
These to terms in $\partial^2$ cancel each other. Thus we get the
associativity equation.
\end{proof}
\end{theorem*}
\subsubsection{{\tt Example}}
Consider the space $\overline{\K_2^2}$. It is a 1-dimensional space. Its
boundary consists from two points, these points are
$\K_{1,1}^2\times\K_2^{1,1}$ and $\K_2^1\times\K_1^2$. We can write, up to a
common sign, $\partial
(\overline{\K_2^2})=\K_2^1\times\K_1^2-\K_{1,1}^2\times\K_2^{1,1}$.
This example explains the signs in these Section.

Now we are ready to give the following definitions:
\begin{definition}
\begin{itemize}
\item[(i)]{\tt A strong homotopy bialgebra structure on a vector space $V$ is a
map of CROCs $\Upsilon\colon\aleph\to\End(V)$,
\item[(ii)]A non-commutative strong homotopy bialgebra structure on a vector
space $V$ is a map of {\it pre}CROCs $\Upsilon\colon\aleph\to\End(V)$.}
\end{itemize}
\end{definition}

\begin{remark}
Note that in the construction above we also define on the associative algebra
$\mho$ a differential, compatible with the product. The differential
comes from the "linear" term in the action of $\partial$ on
$\overline{\aleph}$.
\end{remark}
\section{{\tt The Quillen duality and the (non-)commutative deformations}}
\subsection{{\tt The Quillen duality}}
The classical Quillen duality gives two maps $Q_{C\to
L}\colon\mathrm{Comm}\to\mathrm{Lie}$ from commutative dg algebras to Lie dg
algebras and $Q_{L\to C}\colon\mathrm{Lie}\to\mathrm{Comm}$ from Lie dg
algebras to commutative dg algebras which establish the equivalence of the
derived categories $\mathrm{DComm}$ and $\mathrm{DLie}$. It means, that for
a commutative dg algebra $A^\mb$, the commutative dg algebra $Q_{L\to C}\circ
Q_{C\to L}(A)$ is isomorphic to $A^\mb$ in the derived category, and for
a dg Lie algebra $\g^\mb$, the dg Lie algebra $Q_{C\to L}\circ Q_{L\to
C}(\g^\mb)$ is isomorphic in the derived category to $\g^\mb$.

These functors $Q_{C\to L}$ and $Q_{L\to C}$ are constructed us follows.
For a commutative dg algebra $A^\mb$, consider the free Lie algebra
$\mathrm{Free}((A^\mb [1])^*)$
generated by the dual space $(A^\mb  [1]]^*$. The product in $A^\mb$ is a map
$S^2(A^\mb)\to A^\mb$ where the symmetric square is understood in the graded
sense. Then we have the dual map $\delta\colon (A^\mb)^*\to S^2(A^\mb)^*$.
The map $\delta$ can be considered as a map from the generators of the Lie
algebra $\mathrm{Free}((A^\mb [1])^*)$ to the brackets of the generators. It
turns out from the associativity of the product in $A^\mb$ that the map
$\delta$ can be correctly extended to a differential on $\mathrm{Free}((A^\mb [1])^*)$
of degree $+1$.

On the other hand, for a dg Lie algebra $\g^\mb$, $Q_{L\to C}(\g^\mb)$ is by
definition the chain complex of the Lie algebra $\g^\mb$.

The fact that these to functors define the equivalence of the derived
categories $\mathrm{DComm}$ and $\mathrm{DLie}$ is proven in [Q1,2].

On the other hand, in the same way one can define the functor $Q_{A\to A}$
from associative dg algebras to itself. Namely, for an associative dg
algebra $A^\mb$, consider the free (tensor) associative algebra
$T((A^\mb[1])^*)$ generated by the space $(A^\mb[1])^*$. The dual map to the
product is a map $\delta\colon (A^\mb)^*\to\otimes ^2(A^\mb)^*$, and it
follows from the associativity of the product in $A^\mb$ that $\delta$ can
be continued to a differential in $T((A^\mb[1])^*)$ by the Leibniz rule.
One can prove that $Q^2_{A\to A}(A^\mb)$ is isomorphic to $A^\mb$ in the
derived category.
\begin{remark}
Note that the cohomology of $T((A[1])^*)$ are 0 in all degrees for a degree 0 associative
algebra $A$.
Nevertheless, the double application of this construction gives non-trivial
cohomology. The point is that the corresponding spectral sequence does not
converge to the total cohomology, and we can not use it. Another example of
such a situation: Consider the derivations of $T((A^\mb[1])^*)$. It is
clearly the Hochschild cohomological complex of $A^\mb$, and has non-zero
cohomology for many degree 0 algebras.
\end{remark}
It is clear that the Quillen duallty for associoative algebras is compatible
with the Quillen duality for commutative and for Lie algebras.

\subsection{{\tt Relation with deformation theory}}
In deformation theory, the deformations of an object are described via the
deformation functor. This is a functor on the Artinian algebras constructed
from a dg Lie algebra. This dg Lie algebra is the algebra of the derivations
of the object we deform in a higher sense. In each case, we define this dg
Lie algebra differently. The general prescription is to replace the object
by its resolution and to take the derivations of the resolution. It is a Lie
algebra with the differential equal to the bracket with the differential in
the resolution (which is a distinguished derivation of the resolution).

On the other hand, the 0-th Lie algebra cohomology of this dg Lie algebra
are equal to the commutative algebra of functions on the formal neighborhood
of the object in the moduli space of deformations. To make this claim
rigorous, we should work carefully with the infinite-dimensional objects,
and consider the right completions. For the deformation theory of Riemann
surfaces this claim is proved in [F].

Then, it is clear, that the more right object is the extended dg commutative
algebra, which is the Quillen dual to the deformation dg Lie algebra.
(The cochain complex $Q_{L\to C}(\g^\mb)$ computes the Lie algebra
cohomology).

In Section 2 of the present paper we constructed a dg associative algebra
$\mho$. We can explain its relation with the deformation theory as follows.
Suppose that the {\it extended} (in the sense above) commutative
neighborhood in the moduli space is a part of a bigger non-commutative
space. The there is a map $p\colon A^\mb\to A_0^\mb$ where $A_0^\mb$ is the
extended commutative dg algebra, and $A^\mb$ is the associative dg algebra.
Suppose that $A_0^\mb=A^\mb/[A^\mb,A^\mb]$. The commutant here should be
understood in the sense of the derived functors.

At the moment we do not know how to define the non-commutative algebra
$A^\mb$. We are going to consider this problem in the sequel.
But now we think about the algebra $\mho$ as about the Quillen dual to
$A^\mb$, $\mho=Q_{A\to A}(A^\mb)$. The evidence for this conjecture is that
the product of some elements (the "diagonal" elements) in $\mho$ is looks
very closely to the Maurer-Cartan equation for bialgebras. This informal
conjecture allows us to formulate another one, more rigorous:

\subsubsection*{Conjecture}
{\tt Consider the algebra Quillen dual to $\mho$, $\Omega=Q_{A\to A}(\mho)$.
Consider the quotient $A_0^\mb=\Omega/[\Omega,\Omega]$. Then the 0-th
cohomology of $A_0^\mb$ is isomorphic to the functions on the formal
neighborhood of the initial bialgebra in the moduli space of bialgebras.
Next, the Quillen dual to the dg commutative algebra $A^\mb_0$,
$\g^\mb=Q_{C\to L}(A^\mb_0)$ is the deformation Lie algebra for deformations
of the initial (co)associative bialgebra. It means that the deformation
functor associated with this dg Lie algebra, describes the deformations of
(co)associative bialgebras.}

\begin{remark}
The quotient by the commutant in the Conjecture above can be understood in
the usual sense because the algebra $\Omega$ is free.
\end{remark}
\subsection{{\tt Formality conjectures}}
In the case when the initial bialgebra is $S(V)$ (a free commutative
cocommutative bialgebra), the algebra $\mho$ is quasi-isomorphic to its
cohomology as an associative dg algebra. The corresponding Lie algebra
$\g^\mb$ (constructed in the Conjecture above) in this case is also formal.

In this case, maybe the non-commutative formality (of the algebra $\mho$) is
more simple than the formality of the dg Lie algebra $\g^\mb$.

\comment
\section{{\tt An $L_\infty$-algebra }}

Let $V$ be a vector space. In this Section we define, using the
compactification from Section 1, an $L_\infty$-algebra structure
on a complex {\it quasi-isomorphic } to $\bigoplus_{m,n\ge 1} \Hom(V^{\otimes
n},V^{\otimes m})[-m-n+2]$ (the latter is considered as a complex with 0 differential).
The motivation of this construction is deformation theory of associative
bialgebras, and the Gerstenhaber-Schack complex, which we consider from
our point of view in the next Section.
\subsection{{\tt $L_\infty$ algebras and quasi-$L_\infty$ algebras}}
\begin{definition}
{\tt A quasi-$L_\infty$ algebra structure on a graded vector space $\g^\mb$ is a
{\bf non-homogeneous} vector field $Q$ on $\g^\mb[1]$ such that $[Q,Q]=0$.
The vector field $Q$ is allowed to have both odd and even homogeneous
components}
\end{definition}

Recall that an $L_\infty$ algebra structure on $\g^\mb$ is an odd vector
field $Q$ of degree $+1$ such that $Q^2=0$ (in the odd case,
$Q^2=\frac12[Q,Q]$). In coordinates, such a vector field is a collection of
maps
\begin{equation}\label{eq2.1}
Q_k\colon\wedge^k\g^\mb\to \g^\mb[2-k]
\end{equation}
and the condition $Q^2=0$ becomes the following quadratic relations on
$\{Q_k\}$:
\begin{equation}\label{eq2.2}
\sum_{i+j=k}\Alt_{\alpha_1,\dots,\alpha_k}\pm Q_i(Q_j(\alpha_1,\dots,
\alpha_j),\alpha_{j+1},\dots,\alpha_k)=0
\end{equation}

In the case of quasi-$L_\infty$ algebras, we have maps
\begin{equation}\label{eq2.3}
Q_{k,s}\colon \wedge^k\g^\mb\to\g^\mb[2-k+s],k\ge 1,s\in \mathbb{Z}
\end{equation}
such that the components
\begin{equation}\label{eq2.4}
Q_k:=\sum_{s\in\mathbb{Z}}Q_{k,s}
\end{equation}
satisfy the equation (\ref{eq2.2}).

The numbers $s$ in (\ref{eq2.3}) can be as well positive as negative, and
as well even as odd.

We construct here such a structure on the graded vector space
$\g^\mb=\bigoplus_{m,n\ge 1}\Hom(V^{\otimes n},V^{\otimes m})[-m-n+2]$ where
$V$ is a vector space. In the next Section we show that when $V=A$ is an
associative bialgebra with the product $\Psi\colon A^{\otimes 2}\to A$ and
with the coproducy $\Delta\colon A\to A^{\otimes 2}$ (which are associative
and coassociative, correspondingly, and compatible--see Section 3), the sum
$\Psi +\Delta\in \g^1$ satisfies the generalized Maurer-Cartan equation
$Q(\Psi+\Delta)=0$. Notice that our vector field $Q$ is not lineat, and even
is not polynomial.

\subsection{{\tt The origin of the problem: from strata of codimension 1 to $Q_\sigma$}}

Let
\begin{equation}\label{eq2new.1}
\g^\mb=\bigoplus_{m,n\ge 1}\Hom(V^{\otimes n},V^{\otimes m})[-m-n+2]
\end{equation}
Let $\sigma$ be a configuration of codimension 1 in $\overline{\K(m,n)}$
with the parameters $\ell_1,\ell_2,\{m_i\},\{n_j\}$ (see Section 1.3). Here
we put in a correspondence to $\sigma$ an element
\begin{equation}\label{eq2.5}
Q_\sigma\colon\wedge^{\ell_1+\ell_2}\g^\mb\to\g^\mb[2-\ell_1-\ell_2-(2\ell_1\ell_2-3(\ell_1+\ell_2)+4)]
\end{equation}
The imade of $Q_\sigma$ belongs to $\Hom(V^{\otimes n},V^{\otimes m})$, and
$Q_\sigma$ is non-zero only if its $\ell_1+\ell_2$ arguments are
$$
\Psi_1\in\Hom(V^{\otimes n_1},V^{\otimes \ell_1}),\Psi_2\in\Hom(V^{\otimes
n_2},V^{\otimes \ell_1}),\dots,\Psi_{\ell_2}\in\Hom(V^{\otimes
n_{\ell_2}},V^{\otimes \ell_1});
$$
and
$$
\Theta_1\in\Hom(V^{\otimes \ell_2},V^{\otimes
m_1}),\Theta_2\in\Hom(V^{\otimes \ell_2},V^{\otimes
m_2}),\dots,\Theta_{\ell_1}\in\Hom(V^{\otimes \ell_2},V^{\otimes
m_{\ell_1}}).
$$
That is, we want to define a map $Q_\sigma \in\Hom(V^{\otimes n},V^{\otimes
m})$ starting from
$\Psi_1,\dots,\Psi_{\ell_2};\Theta_1,\dots,\Theta_{\ell_1}$ as above. The
construction is as follows:

First define
\begin{multline}\label{eq2.6}
F_\sigma(v_1\otimes\dots\otimes v_n)\colon =
\Psi_1(v_1\otimes\dots\otimes
v_{n_1})\bigotimes\Psi_2(v_{n_1+1}\otimes\dots\otimes
v_{n_2})\bigotimes\dots\\ \bigotimes
\Psi_{\ell_2}(v_{n_{\ell_2-1}+1}\otimes\dots\otimes v_{n_{\ell_2}})
\in V^{\otimes \ell_1\ell_2}
\end{multline}
Now we apply $\{\Theta_i\}$ to this element in $V^{\otimes \ell_1\ell_2}$:
we define an element $G_\sigma\colon V^{\otimes \ell_1\ell_2}\to V^{\otimes m}$ as follows:
\begin{multline}\label{eq2.7}
G_\sigma(v_1\otimes v_2\otimes\dots\otimes v_{\ell_1\ell_2}):=
\Theta_1(v_1\otimes v_{\ell_1+1}\otimes\dots\otimes
v_{\ell_1(\ell_2-1)+1})\bigotimes\\
\Theta_2(v_2\otimes v_{\ell_1+2}\otimes\dots\otimes
v_{\ell_1(\ell_2-1)+2})\bigotimes\dots
\bigotimes\Theta_{\ell_1}(v_{\ell_1}\otimes v_{2\ell_1}\otimes\dots\otimes
v_{\ell_1\ell_2})\in V^{\otimes m}.
\end{multline}
Define now
\begin{equation}\label{eq2.8}
Q_\sigma(v_1\otimes \dots\otimes v_n):=G_\sigma\circ
F_\sigma(v_1\otimes\dots\otimes v_n)
\end{equation}
We denote by the same symbol $Q_\sigma$ the image of $Q_\sigma$ after the
alternation of $\{v_1,\dots, v_n\}$.

Now let us compute the degrees. We have:
$$
\deg \Psi_i=n_i+\ell_1-2,\ \deg\Theta_i=m_i+\ell_2-2,\ \deg Q_\sigma=m+n-2
$$
We are interesting in the "defect"
$$
N_\sigma:=\deg Q_\sigma-(\sum_{i=1}^{\ell_2}\deg \Psi_i+\sum_{i=1}^{\ell_1}
\deg \Theta_i+(2-\ell_1-\ell_2))
$$
Recall that an $L_\infty$-algebra structure on a graded vector space
$\g^\mb$ is an odd vector
field $Q$ of degree $+1$ such that $Q^2=0$ (in the odd case,
$Q^2=\frac12[Q,Q]$). In coordinates, such a vector field is a collection of
maps
\begin{equation}\label{eq2.1}
Q_k\colon\wedge^k\g^\mb\to \g^\mb[2-k]
\end{equation}
and the condition $Q^2=0$ becomes the following quadratic relations on
$\{Q_k\}$:
\begin{equation}\label{eq2.2}
\sum_{i+j=k}\Alt_{\alpha_1,\dots,\alpha_k}\pm Q_i(Q_j(\alpha_1,\dots,
\alpha_j),\alpha_{j+1},\dots,\alpha_k)=0
\end{equation}

If $N_\sigma$ would be equal to 0, we would have the {\it pure} shift of
degrees (as in the case of $L_\infty$ algebras). It is clear that
\begin{equation}\label{eq2.9}
N_\sigma=-(2\ell_1\ell_2-3(\ell_1+\ell_2)+4)
\end{equation}

This number $N_\sigma$ is equal to 0 only in the following two cases:
\begin{lemma}
For integral $\ell_1,\ell_2$, $2\ell_1\ell_2-3(\ell_1+\ell_2)+4=0$ iff
$(\ell_1,\ell_2)=(1,1)$ or $(\ell_1,\ell_2)=(2,2)$.
\qed
\end{lemma}

We see, therefore, that we can not consider the maps $Q_\sigma$ as
the components of an odd vector field $Q$ of degree $+1$ on $\g^\mb [1]$.
Nevertheless, we want to modify the construction, and to deduce the
cohomological equation $Q^2=0$ from the equation $\partial^2=0$
where $\partial$ is the boundary operator in the stratified space
$\overline{\K(m,n)}$. For this we replace the graded space $\g^\mb$ by a
quasi-isomorphic complex and define an $L_\infty$-algebra structure on this
extended complex. It turns out, however, that the induced Lie algebra
structure on the cohomology (isomorphic as a graded vector space to to $\g^\mb$)
 is {\it not} equivalent as an
$L_\infty$-algebra to the $L_\infty$-structure on the extended complex.

\subsection{{\tt The Construction}}

Denote by $C^\mb_{m,n}$ the cell copmplex of the stratified space
$\overline{\K(m,n)}$. By definition, $C^{-k}_{m,n}$ is the vector space over
$\mathbb{R}$ generated formally by the strata of dimension k. The
differential $\partial\colon C_{m,n}^i\to C_{m,n}^{i+1}$ is defined in the
natural way. This is a differential of degree +1. It is clear that
$\X^i(C_{m,n}^\mb)\ne 0$ only for $i=0$, and $\X^0(C_{m,n}^\mb)=\mathbb{R}$.
It follows from the fact that the space $\overline{\K(m,n)}$ is
contractible.

Now let $V$ be a vector space. Consider instead of
$$
\g^\mb=\bigoplus_{m,n\ge 1}\Hom(V^{\otimes n},V^{\otimes m})[-m-n+2]
$$
the complex
\begin {equation}\label{eq2new10}
\wtilde{\g}^\mb=\bigoplus_{m,n\ge 1}C^\mb_{m,n}\otimes_{\mathbb{R}}\Hom(V^{\otimes n},V^{\otimes m})[-m-n+2]
\end{equation}
It is clear that $\X^\mb(\wtilde{\g}^\mb)\simeq\g^\mb$.

We are going to construct an $L_\infty$-algebra structure on
$\wtilde{\g}^\mb$ using the operations $Q_\sigma$ ($\sigma$ is a stratum of
codimension 1) defined in the previous Subsection, and some their
generalizations.

Denote the complex $C_{m,n}^\mb\otimes_{\mathbb{R}}\Hom(V^{\otimes
n},V^{\otimes m})[-m-n+2]$ by $\wtilde{\Hom}^\mb(V^{\otimes n},V^{\otimes m})$.
We have: $\X^\mb \wtilde{\Hom}^\mb(V^{\otimes n},V^{\otimes m})\simeq
\Hom(V^{\otimes n},V^{\otimes m})[-m-n+2]$. As
$\dim\overline{\K(m,n)}=m+n-3$, the complex $\wtilde{\Hom}^\mb(V^{\otimes n},V^{\otimes
m})$ has degrees from 1 to $m+n-2$.

The elements of degree 1 in $\wtilde{\Hom}^\mb(V^{\otimes n},V^{\otimes m})$
are corresponded
to the strata of codimension 0 in $\overline{\K(m,n)}$ (there is only 1 such
stratum), the elements of degree 2 in $\wtilde{\Hom}^\mb(V^{\otimes n},V^{\otimes
m})$ are corresponded to the strata of codimension 1 in
$\overline{\K(m,n)}$, and so on.

Now we attach to each stratum $\sigma$ of codimension 1 in
$\overline{\K(m,n)}$ with the parameters $\ell_1,\ell_2,\{m_i\},\{n_j\}$
(see Section 1.3) an operation
\begin{multline}\label{eq2new11}
\wtilde{Q_\sigma}\colon \wtilde{\Hom}^1(V^{\otimes n_1},V^{\otimes \ell_1})
\otimes\dots\otimes\wtilde{\Hom}^1(V^{\otimes n_{\ell_2}},V^{\otimes
\ell_1})\otimes\\
\otimes \wtilde{\Hom}^1(V^{\otimes \ell_2},V^{\otimes
m_1})\otimes\dots\otimes\wtilde{\Hom}^1(V^{\otimes \ell_2},V^{\otimes m_{\ell_1}})\longrightarrow
\wtilde{\Hom}^2(V^{\otimes n},V^{\otimes m})[2-\ell_1-\ell_2]
\end{multline}
$(m=\sum_{i=1}^{\ell_1}m_i,\ \ n=\sum_{j=1}^{\ell_2}n_j)$.

We define the map $\wtilde{Q_\sigma}$ exactly by the formulas
(\ref{eq2.6}), (\ref{eq2.7}), (\ref{eq2.8}) where
all $\Psi_i$ and $\Theta_j$ have the grading (degree) +1 in the
$\wtilde{\Hom}$'s complexes. We define the result
$$
\wtilde{Q_\sigma}(\Psi_1\otimes\dots\otimes\Psi_{\ell_2};\Theta_1\otimes\dots\otimes\Theta_{\ell_1})
$$
as the element in
$\wtilde{\Hom}^2(V^{\otimes n},V^{\otimes m})$ of the grading 2
(corresponded to the stratum $\sigma$ of codimension 1 in
$\overline{\K(m,n)}$) equal to
\begin{multline}\label{eq2new4}
\wtilde{Q_\sigma}(\Psi_1\otimes\dots\otimes\Psi_{\ell_2};\Theta_1\otimes\dots\otimes\Theta_{\ell_1})=\\
=\sigma\bigotimes Q_\sigma (\Psi_1\otimes\dots\otimes\Psi_{\ell_2};\Theta_1\otimes\dots\otimes\Theta_{\ell_1})
\end{multline}

The advantage of this definition comparably with our definition of
$Q_\sigma$ in Section 2.1 is that $\wtilde{Q_\sigma}$has the correct (in the
sense of $L_\infty$-algebras)degree. Namely, each $\Psi_i$ and $\Theta_j$
have degree +1, whence
$\wtilde{Q_\sigma}(\Psi_1\otimes\dots\otimes\Psi_{\ell_2};\Theta_1\otimes\dots\otimes\Theta_{\ell_1})$
has degree +2. Therefore,
\begin{multline}\label{eq2new5}
\deg \wtilde{Q_\sigma}(\Psi_1\otimes\dots\otimes\Psi_{\ell_2};\Theta_1\otimes\dots\otimes\Theta_{\ell_1})
=\sum_{i=1}^{\ell_2}\deg \Psi_i+\sum_{j=1}^{\ell_1}\deg
\Theta_j+(2-\ell_1-\ell_2)
\end{multline}

Now we want to extend this definition of $\wtilde{Q_\sigma}$to the elements
of arbitrary degree in the $\wtilde{\Hom}$'s complexes. For this we
introduce the technique we call {\it the glueing } in the next Subsection.

\subsection{{\tt The Glueing}}
To proceed in this way, we need to know the answer on the following
question:
\begin{question}
{\tt Does there exist a natural map which attaches to
a stratum of codimension $k_1$ in $\overline{\K(\ell_1,n_1)}$,..., a stratum
of codimension $k_{\ell_2}$ in $\overline{\K(\ell_1,n_{\ell_2})}$, and a
stratum of codimension $p_1$ in $\overline{\K(m_1,\ell_2)}$,..., a stratum
of codimension $p_{\ell_1}$ in $\overline{\K(m_{\ell_1},\ell_2)}$, a stratum
of codimension $(\sum_{i=1}^{\ell_2}k_i+\sum_{j=1}^{\ell_1}p_j)+1$ in
$\overline{\K(\sum_{j=1}^{\ell_1}m_j,\sum_{i=1}^{\ell_2}n_i)}$? }
\end{question}
If such a way exists, we could directly generalize the construction of
$\wtilde{Q_\sigma}$ from the previous Subsection, and to define a map
\begin{multline}\label{eq2new5}
\wtilde{Q_\sigma}\colon \wtilde{\Hom}^{k_1+1}(V^{\otimes n_1},V^{\otimes \ell_1})
\otimes\dots\otimes\wtilde{\Hom}^{k_{\ell_2}+1}(V^{\otimes n_{\ell_2}},V^{\otimes
\ell_1})\otimes\\
\otimes \wtilde{\Hom}^{p_1+1}(V^{\otimes \ell_2},V^{\otimes
m_1})\otimes\dots\otimes\wtilde{\Hom}^{p_{\ell_1}+1}(V^{\otimes \ell_2},V^{\otimes
m_{\ell_1}})\longrightarrow\\
\longrightarrow
\wtilde{\Hom}^{\sum_{i=1}^{\ell_2} k_i+\sum_{j=1}^{\ell_1} p_j +2}(V^{\otimes n},V^{\otimes m})[2-\ell_1-\ell_2]
\end{multline}
$(m=\sum_{i=1}^{\ell_1}m_i,\ \ n=\sum_{j=1}^{\ell_2}n_j)$.

To find a solution to the Question above,
consider the natural map
\begin{multline}\label{eq2new6}
\phi_\sigma\colon\K_{\ell_2}^{m_1,\dots,m_{\ell_1}}\times
\K_{n_1,\dots,n_{\ell_2}}^{\ell_1}\longrightarrow\\
\K(m_1,\ell_2)\times\dots\times\K(m_{\ell_1},\ell_2)
\times\K(\ell_1,n_1)\times\dots\times\K(\ell_1,n_{\ell_2})
\end{multline}
(Here $\sigma$ is a stratum of codimension 1 in $\overline{\K(m,n)}$ with
the parameters $\ell_1,\ell_2,\{m_i\},\{n_j\}$).
\begin{lemma}
\begin{itemize}
\item[(i)] The map $\phi_\sigma$ can be canonically extended to a map
\begin{multline}\label{eq2new7}
\overline{\phi_\sigma}\colon\overline{\K_{\ell_2}^{m_1,\dots,m_{\ell_1}}\times
\K_{n_1,\dots,n_{\ell_2}}^{\ell_1}}\longrightarrow\\
\overline{\K(m_1,\ell_2)}\times\dots\times\overline{\K(m_{\ell_1},\ell_2)}
\times\overline{\K(\ell_1,n_1)}\times\dots\times\overline{\K(\ell_1,n_{\ell_2})}
\end{multline}
where $\overline{\K_{\ell_2}^{m_1,\dots,m_{\ell_1}}\times
\K_{n_1,\dots,n_{\ell_2}}^{\ell_1}}$ denotes the closure of the stratum of
codimension 1 $\K_{\ell_2}^{m_1,\dots,m_{\ell_1}}\times
\K_{n_1,\dots,n_{\ell_2}}^{\ell_1}$ in
$\overline{\K_{n_1,\dots,n_{\ell_2}}^{m_1,\dots,m_{\ell_1}}}$.
\item[(ii)] For each stratum (of higher codimension) $\omega\in\overline{\K_{\ell_2}^{m_1,\dots,m_{\ell_1}}\times
\K_{n_1,\dots,n_{\ell_2}}^{\ell_1}}$, the projection
$p_{\alpha}(\overline{\phi_\sigma}(\omega))$ to the $\alpha$'s factor in
$\overline{\K(m_1,\ell_2)}\times\dots\times\overline{\K(m_{\ell_1},\ell_2)}
\times\overline{\K(\ell_1,n_1)}\times\dots\times\overline{\K(\ell_1,n_{\ell_2})}$
is the only one stratum (in the stratification described in Section 1) of
some codimension $d_{\alpha}$
\item[(iii)] The codimension $d_{\omega}$ of $\omega$ in
$\overline{\K_{\ell_2}^{m_1,\dots,m_{\ell_1}}\times
\K_{n_1,\dots,n_{\ell_2}}^{\ell_1}}$ is equal to the sum
\begin{equation}\label{eq2new8}
d_{\omega}=\sum_{all\ factors\ \alpha}d_{\alpha}
\end{equation}
where in the sum $\alpha$ runs trough all factors of $\overline{\K(m_1,\ell_2)}\times\dots\times\overline{\K(m_{\ell_1},\ell_2)}
\times\overline{\K(\ell_1,n_1)}\times\dots\times\overline{\K(\ell_1,n_{\ell_2})}$.
\end{itemize}
\end{lemma}

\endcomment
\subsection*{Acknowledgements}
The author is grateful to Borya Feigin and to Maxim Kontsevich for many
discussions. Discussions with Andrey Losev on the CROCs of associative
algebras and on Stasheff associahedrons were very useful for me. I am
grateful very much to Giovanni Felder for organizing and supervising my postdoc
at the ETH (Zurich). I would like to express my gratitude to the ETH and to
the IHES (Bures-sur-Yvette) where a part of this work was done for the
hospitality and the excellent working conditions.

\bigskip
\bigskip
Dept. of Math., ETH-Zentrum, 8092 Zurich, SWITZERLAND\\
e-mail: {\tt borya@mccme.ru, borya@math.ethz.ch}


\begin{thebibliography}{999}
\bibitem[F]{F}
B.~Feigin, The conformal field theory from the view of the cohomology theory
of Lie algebras, ICM-90 (Kyoto) , Math. Society of Japan, Tokyo, 1990
\bibitem[K1]{K1}
M.~Kontsevich, Deformation quantization of Poisson manifolds I, preprint
q-alg/9709040
\bibitem[K2]{K2}
M.~Kontsevich, Operads and Motives in Deformation Quantization, preprint
math.QA/9904055
\bibitem[M]{M}
M.~Markl, A resolution (minimal model) of the PROP for bialgebras, preprint
math.AT/0209007
\bibitem[Q1]{Q1}
D.G.~Quillen, Homotopical algebra, Lect. Notes in Math. 43, Springer-Verlag,
1967
\bibitem[Q2]{Q2}
D.G.~Quillen, Rational Homotopy Theory, Ann. of Math., vol.90, 1969
\bibitem[T]{T}
D.~Tamarkin, Another proof of M.Kontsevich formality theorem, preprint
math.QA/9803025
\end{thebibliography}
\end{document}